\pdfoutput=1
\RequirePackage{ifpdf}
\ifpdf 
\documentclass[pdftex]{sigma}
\else
\documentclass{sigma}
\fi

\def\hti{\widetilde{\mathfrak{h}}}
\def\s{\mathfrak{s}}
\def\n{\mathfrak{n}}
\def\h{\mathfrak{h}}
\def\Z{\mathbb Z}

\def\ga{\hat{\g}}
\def\g{\mathfrak{g}}
\def\gt{\widetilde{\g}}
\def\Span{\text{span}}
\def\N{\mathbb N}

\usepackage{mathtools} 

\usepackage{array} 
\usepackage{booktabs} 

\numberwithin{equation}{section}

\newtheorem{Theorem}{Theorem}[section]
\newtheorem*{Theorem*}{Theorem}

\newtheorem{Proposition}[Theorem]{Proposition}

\theoremstyle{definition}

\newtheorem{Remark}[Theorem]{Remark}

\begin{document}

\allowdisplaybreaks

\newcommand{\arXivNumber}{2511.12284}

\renewcommand{\thefootnote}{}

\renewcommand{\PaperNumber}{055}

\FirstPageHeading

\ShortArticleName{Leading Terms of Relations on a Level 5 Module over the Twisted Affine Lie Algebra \smash{$A_2^{(2)}$}}

\ArticleName{Leading Terms of Relations on a Level 5 Module\\ over the Twisted Affine Lie Algebra $\boldsymbol{A_2^{(2)}}$\footnote{This paper is a~contribution to the Special Issue on Recent Advances in Vertex Operator Algebras in honor of James Lepowsky. The~full collection is available at \href{https://sigma-journal.com/Lepowsky.html}{https://sigma-journal.com/Lepowsky.html}}}

\Author{Stefano CAPPARELLI~$^{\rm a}$, Arne MEURMAN~$^{\rm b}$ and Mirko PRIMC~$^{\rm c}$}

\AuthorNameForHeading{S.~Capparelli, A.~Meurman and M.~Primc}

\Address{$^{\rm a)}$~Dipartimento SBAI, Universit\`a La Sapienza, Rome, Italy}
\EmailD{\mail{stefano.capparelli@uniroma1.it}}

\Address{$^{\rm b)}$~Department of Mathematics, Lund University, Lund, Sweden}
\EmailD{\mail{arne.meurman@math.lu.se}}

\Address{$^{\rm c)}$~Faculty of Science, University of Zagreb, Zagreb, Croatia}
\EmailD{\mail{primc@math.hr}}

\ArticleDates{Received November 18, 2025, in final form May 18, 2026; Published online June 02, 2026}

\Abstract{One of the starting points of this work was the duality of Borcea relating standard~level~$k$ representations of \smash{$A_1^{(1)}$} and level $2k+1$ of \smash{$A_2^{(2)}$}. For $k=1$, the combinatorial bases in both cases yield the two Capparelli identities and we wanted to see if there is a~correspondence between the bases in terms of partitions for all $k\in\mathbb N$. By using the vertex operator relations in the principal picture for level $5$ standard \smash{$A_2^{(2)}$}-modules, we~reduce a~spanning set of Poincar\'e--Birkhoff--Witt-type vectors in $L(5\Lambda_0)$ by removing the leading terms of relations and rendering a list of 34 ``difference'' conditions for partitions. Using computer programs, we enumerated the partitions satisfying these conditions and obtained a~truncated generating series agreeing with the principally specialized character for all powers of $q$ up to $41$. Although our list of leading terms is incomplete, our results show that the corresponding combinatorial identity for \smash{$L_{A_2^{(2)}}(5\Lambda_0)$} drastically differs from the one for the Borcea dual \smash{$L_{A_1^{(1)}}(2\Lambda_0)$}.}


\Classification{17B69; 11P84}

\begin{flushright}
\begin{minipage}{85mm}
\it The authors are deeply grateful to Jim Lepowsky\\ for long-term collaboration,
support, and friendship
\end{minipage}
\end{flushright}

\renewcommand{\thefootnote}{\arabic{footnote}}
\setcounter{footnote}{0}

\section{Introduction}

In \cite{LW82}, J.~Lepowsky and R.L.~Wilson discovered that the classical Rogers--Ramanujan identities are connected
to the level 3 standard modules over the affine Lie algebra of type \smash{$A_1^{(1)}$}. They introduced a
Heisenberg subalgebra $\mathfrak{s}$ and ``$Z$-operators'' commuting with $\mathfrak{s}$, hence acting on
the vacuum space of $\mathfrak{s}$. The $Z$-operators were shown to satisfy quadratic relations, and this
led to a~combinatorial description of a basis of the vacuum space in terms of the partitions from the sum
sides of the Rogers--Ramanujan identities. The product side in \cite{LM78} resulted from Lepowsky's
numerator formula combined with the principal specialization of the Weyl--Kac character formula. In
\cite{MP87} two of the authors, A.~Meurman and M.~Primc, obtained Poincar\'e--Birkhoff--Witt type bases for all
of the higher level standard modules over \smash{$A_1^{(1)}$}. This gave Lie theoretic proofs of the combinatorial
identities of B.~Gordon, G.E.~Andrews, and D.~Bressoud~\mbox{\cite{A67a, A67b, Br79, Br80, G61}}.

 Other connections between Lie algebra
representations and combinatorial identities have been obtained in, e.g., \cite{CMPP21, DK21,FKLMM99, KR23, MP99, PT24, PS16,TT20}.

In \cite{C93}, one of the authors, S.~Capparelli, studied the level 3 modules over the Lie algebra of type \smash{$A_2^{(2)}$} and discovered by the Lepowsky--Wilson method combinatorial identities of Rogers--Ramanujan type. These were the first new such identities discovered via Lie algebra representation theory. They were subsequently proved with combinatorial methods in~\cite{A92}, and~with~Lie-algebraic methods in~\cite{C96, TX95}. In~\cite{N14}, D.~Nandi studied the level~4 modules and conjectured three new Rogers--Ramanujan type identities. These were later proved in~\cite{TT19}. In~this paper, we study in a similar way the level 5 module $L(5\Lambda_0)$ over \smash{$A_2^{(2)}$}. Using the methods of Capparelli and Nandi, we obtain in Section~\ref{S6} a list of 34 conditions on the ``sum'' side of the associated combinatorial identity. For partitions of size at most~$48$, we have equality between the sum and product sides, except that at partitions of 42 respectively 48 we lack one partition each in our conditions. Thus our list is not yet complete, but we would still like to publish our results to this point for the benefit of other researchers.

Our study was motivated in part by a duality observed by J.~Borcea. In~\cite{Bo02}, he noticed a~correspondence between the characters of the level $k$ modules over the Lie algebra \smash{$A_1^{(1)}$} and the level $(2k+1)$-modules over \smash{$A_2^{(2)}$}. In Section \ref{S7}, we consider the $k=2$ case of Borcea's duality, and present the basis and combinatorial identity associated to the \smash{$A_1^{(1)}$}-module $L(2\Lambda_0)$. Our results show that the sum sides associated to the two modules differ drastically.

The paper is organized as follows. In Section~\ref{S2}, we recall the basic set up of the vertex operator construction of the basic module of the principally graded realization of the affine algebra \smash{$A_2^{(2)}$}. In particular, we used the construction given in \cite{KKLW81}. This is a special case of the general theory of vertex operator calculus as in \cite{DL93,FLM88, L85}. We also present the affine algebra $\hat{\mathfrak{g}}$ of type \smash{$A_2^{(2)}$} with generators and relations associated to a generalized Cartan matrix as in~\cite{K90}.
In Section~\ref{S3}, we recall the basic notions of standard modules for the algebra $\hat{\mathfrak{g}}$. There is a unique, up to equivalence, standard module of level 1 which we denote $U$ which has highest weight $\Lambda_0$. Our focus will be on a level~5 module with highest weight $5\Lambda_0$ which we view as submodule of the tensor product of copies of $U$. We also use Lepowsky's numerator formula and the principal specialization of the Weyl--Kac character formula.

In Section~\ref{S4}, we present the cubic relations Theorems~\ref{Rrelation} and~\ref{Srelation} and the higher order relations Propositions~\ref{Sright}--\ref{Sleftright} that we need to obtain the leading terms in Section \ref{S6}. Section \ref{S5} is devoted to a detailed explanation of the proof of one leading term $X(-(5,5,2,2))v_\Lambda$. This results from a linear combination of 8 relations acting on the highest weight vector $v_\Lambda$. We use \textsc{Maple} to verify the existence of the linear combination, see \href{http://arxiv.org/abs/2511.12284}{arXiv:2511.12284}.

In Section~\ref{S6}, we formulate our main result about leading terms for the module $L(5\Lambda_0)$. This is given by a list of 34 conditions that the partition $\mu$ in a generating vector $\alpha(-\lambda)X(-\mu)v_\Lambda$ has to satisfy in order to remain in the spanning set. As in Section~\ref{S5}, the linear relations are obtained by \textsc{Maple} calculations, and the \textsc{Maple} worksheets used are available for download. Finally, Section~\ref{S7} is devoted to the case $k=2$ of the duality of Borcea.

\section{Preliminaries}\label{S2}

\subsection[The algebra k]{The algebra $\boldsymbol{\mathfrak{k}}$}
Let $\Phi$ be the $A_2$ root system with basis $\Delta=\{\alpha_1,\alpha_2\}$. Let $L = \mathbb{Z}\alpha_1\oplus\mathbb{Z}\alpha_2$ be the root lattice of type $A_2$ equipped with a symmetric $\mathbb{Z}$-bilinear form such that $\langle \alpha_i,\alpha_i\rangle=2$, $i=1,2$, and
$\langle \alpha_1,\alpha_2\rangle=-1$.
Let $\nu$ be the automorphism of $L$ of order $6$ acting as a rotation on the root system by \smash{$\frac{\pi}{3}$}. On the basis elements of $\Delta$, $\nu(\alpha_1)=\alpha_1+\alpha_2$, $\nu(\alpha_2)=-\alpha_1$. Note that
\begin{equation*}
\sum_{p\in \mathbb{Z}_6}\nu^p\alpha=0 \qquad \text{for all } \alpha\in L.
\end{equation*}

{\samepage Let $\mathfrak{h}=\mathbb{C}\otimes _{\mathbb{Z}}L$. Linearly extend the form to all of $\mathfrak{h}$. Let $\omega$ be a primitive 6th root of 1, and for concreteness, choose
$\omega={\rm e}^{{\rm i}\pi/3}=\frac 12+\frac{\sqrt{3}}{2}{\rm i}$.
For $n\in \mathbb{Z}$, set
\begin{equation*}
\mathfrak{h}_{(n)}=\{x\in \mathfrak{h}\mid\nu x=\omega^n x\}
\end{equation*}}%
and we have
\begin{equation*}
\mathfrak{h}=\coprod_{p\in \mathbb{Z}_6}\mathfrak{h}_{(p)}.
\end{equation*}
Note that $\mathfrak{h}_{(n)}=0$ unless $n\equiv \pm1 \pmod{6}$.

Viewing $\mathfrak{h}$ as an abelian Lie algebra, construct the $\nu$-twisted affine Lie algebra
\[
\hti=\hti[\nu]=\coprod_{n\in\mathbb{Z}} \bigl(\mathfrak{h}_{(n)}\otimes t^{n/6}\bigr)\oplus \mathbb{C} c\oplus \mathbb{C} d
\]
with brackets
\begin{gather*}
\bigl[x\otimes t^{i/6},y\otimes t^{j/6}\bigr]
 = \frac{i}{6}\langle x,y\rangle \delta_{i+j,0}c,\qquad
\bigl[d,x\otimes t^{i/6}\bigr]
 = i x\otimes t^{i/6},\\
\bigl[c,x\otimes t^{i/6}\bigr]
 = [c,d]=0.
\end{gather*}

for all $i,j,\in\mathbb{Z}$, $x\in \mathfrak{h}_{(i)}$, $y\in \mathfrak{h}_{(j)}$.
Consider the commutator subalgebra
\begin{equation*}
 \s=\coprod_{n\in\mathbb{Z}, n\ne 0}\bigl(\mathfrak{h}_{(n)} \otimes t^{n/6}\bigr)\oplus \mathbb{C} c
\end{equation*}
(the Heisenberg subalgebra)
and the subalgebras
\begin{gather*}
 \s_{\pm}=\coprod_{n\in\mathbb{Z},\, \pm n> 0}\bigl(\mathfrak{h}_{(n)} \otimes t^{n/6}\bigr)\oplus \mathbb{C} c,
\qquad
 \mathfrak{b}=\s_+\oplus \mathbb{C} c\oplus \mathbb{C} d.
\end{gather*}

Consider $\mathbb{C}$ as a 1-dimensional $\mathfrak{b}$-module on which $\s_+$ and $d$ act trivially and $c$ acts as the identity. Form the induced module
\[
S=\mathcal{U}\bigl(\hti[\nu]\bigr)\otimes_{\mathcal{U}(\mathfrak{b})}\mathbb{C}\simeq \mathcal{S}(\s_-),
\]
where $\mathcal{S}(\s_-)$ is the symmetric algebra on $\s_-$. This is an irreducible module for the Heisenberg subalgebra.

The action of $d$ defines a $\mathbb{Z}$-grading on $S$:
\[
S=\coprod_{n\in -\mathbb{N}}S_n.
\]
For $\alpha\in \mathfrak{h}$ and $n\in \mathbb{Z}$, define \smash{$\alpha_{(n)}$} as the projection of $\alpha$ onto \smash{$\mathfrak{h}_{(n)}$}. Then \smash{$\alpha_{(n)}=0 $}
unless $n\equiv \pm1 \pmod 6$. For $n\in \mathbb{Z}$, define the operator $\alpha(n)=\alpha_{(n)}\otimes t^{n/6}$ on $S$.

In the generating series below, we shall use commuting formal variables $w, w_1,w_2,\dots$. These are related to the formal variables $z$, $z_j$ in \cite{C92,N14} by \smash{$w=z^{1/6}$}, \smash{$w_j = z_j^{1/6}$}.
	
	Set
	\begin{equation*}
		F = \{j\in\mathbb{Z}\mid j\equiv \pm 1 \ (\bmod \ 6)\},
	\end{equation*}
	the set of indices in the Heisenberg generators $\alpha(n)$. For $i\in\mathbb{Z}$, $\alpha\in\Phi$, we shall need the generating series
	\begin{equation*}
		E^{\pm}\bigl(\nu^i\alpha;w\bigr) = \exp\Biggl(6\sum_{n\in F,\pm n > 0}\alpha(n)\frac{\bigl(\omega^{-i}w\bigr)^{-n}}{n}\Biggr).
	\end{equation*}
	Then the vertex operator construction of the level 1 standard module \cite{KKLW81} is given by{\samepage
	\begin{equation}\label{vertexOperatorFormula}
		X(\alpha;w) = KE^-(-\alpha;w)E^+(-\alpha;w),
	\end{equation}
	where $K = (1+\omega)/36$.}

We have the commutation relation, for $\alpha,\beta\in \h$,
\begin{equation*}
 E^+(\alpha;w_1)E^-(\beta;w_2)=E^-(\beta;w_2)E^+(\alpha;w_1)\prod_{p\in\Z_6}\biggl(1-\omega^{-p}\frac{w_2}{w_1} \biggr)^{\langle \nu^p\alpha,\beta \rangle}.
\end{equation*}

We have
\begin{gather*}
X(\nu\alpha;w)=\lim_{w\to \omega^{-1}w}X(\alpha;w),
\qquad
DX(\alpha;w)=-[d,X(\alpha;w)],
\end{gather*}
where $D=w\frac{{\rm d}}{{\rm d}w}$.

 For $\alpha\in L$, define
\[
\alpha(w)=\sum_{n\in \mathbb{Z}}\alpha(n)w^{-n-6}.
\]

Define the coefficients $X(\alpha;n)$ as follows:
\[
X(\alpha;w)=\sum_{n\in\mathbb{Z}}X(\alpha;n)w^{-n}.
\]
Then $X(\alpha;n)$ is a well-defined operator on $S$ of degree $n$.

Notice that with this notation
\[
X\bigl(\nu^k\alpha;n\bigr)=\omega^{kn}X(\alpha;n)
\]
for all $n,k\in\mathbb{Z}$.
The most important commutator formula for us is when $\alpha\in \Phi$, i.e., when $\langle\alpha,\alpha\rangle=2$. For this,
\begin{align*}
[X(\alpha;w_1),X(\alpha;w_2)]
={}&\frac{\omega^2}{6}X(\nu\alpha;w_2)\delta\bigl(\omega^{-2}x\bigr)
 -\frac{\omega^2}{6}X\bigl(\nu^{-1}\alpha;w_2\bigr)\delta\bigl(\omega^{2}x\bigr)
\\
&{}
{+}\,\frac{\omega}{36}cD\delta(-x)
-\frac{\omega}{6}w_2^{6}\alpha(w_2)\delta(-x),
\end{align*}
where \smash{$x=\frac{w_2}{w_1}$}.

From now on, we fix $\alpha = \alpha_1$. Then the general element in $\Phi$ can be denoted $\nu^i\alpha$ for some $i\in\mathbb{Z}_6$.
Writing sometimes $X(n)=X(\alpha;n)$, we have in particular
\begin{gather}
[\alpha(m),\alpha(n)] = \frac{m}{6}\,\delta_{m+n,0}\,c
\qquad \text{if } m,n\equiv \pm 1 \pmod 6, \label{eq:comm-aa}\\
[\alpha(m), X(n)] = X(m+n)
\qquad \text{if } m\equiv \pm 1 \pmod 6, \label{eq:comm-ax}\\
[X(m),X(n)] = \frac{\omega^{2}}{6}\bigl(\omega^{\,n-m}-\omega^{\,m-n}\bigr)X(m+n) \notag\\
\hphantom{[X(m),X(n)] =}{}\,
- \frac{\omega}{6}(-1)^{m} \alpha(m+n)
+ \delta_{m+n,0}\,\frac{\omega}{36}(-1)^{m} m\,c. \label{eq:comm-xx}
\end{gather}

Notice
\begin{equation*}
 \omega^{n-m}-\omega^{m-n}=\begin{cases}
 0 & \text{if } n-m\equiv 0,3 \pmod 6, \\
 \sqrt{3} {\rm i} & \text{if } n-m\equiv 1,2 \pmod 6, \\
 -\sqrt{3}{\rm i} & \text{if } n-m\equiv 4,5 \pmod 6.
\end{cases}
\end{equation*}

The commutation relations \eqref{eq:comm-aa}--\eqref{eq:comm-xx} show (together with $[\mathfrak{k},c]=0$) that
\begin{equation*}
	\mathfrak{k} = \operatorname{span}\{X(n),\,\alpha(j),\,c\mid n\in\mathbb{Z},\, j\in F\}
\end{equation*}
forms a Lie algebra. We shall next present $\mathfrak{k}$ as a Kac--Moody algebra of type \smash{$A_2^{(2)}$}.

\subsection[The algebra A\_2\^{}(2)]{The algebra $\boldsymbol{A_2^{(2)}}$}
The Lie algebra of type \smash{$A_2^{(2)}$} is $\hat{\mathfrak{g}}$ associated to the
generalized Cartan matrix
\[
A=(a_{ij})=\begin{pmatrix}
 \hphantom{-} 2 & -4 \\
 -1 & \hphantom{-} 2
 \end{pmatrix}
\]
and generators $h_0$, $h_1$, $e_0$, $e_1$, $f_0$, $f_1$ with relations
\begin{gather*}
[h_i,h_j]= 0,\qquad
[h_i,e_j]= a_{ij}e_j,\qquad
[h_i,f_j]= -a_{ij}f_j,\qquad
[e_i,f_j]= \delta_{ij}h_i,\\
(\operatorname{ad} e_i)^{-a_{ij}+1}e_j
= 0 \qquad \text{for } i\ne j,\qquad
(\operatorname{ad} f_i)^{-a_{ij}+1}f_j
= 0 \qquad \text{for } i\ne j,
\end{gather*}
where $i,j\in \{0,1\}$ and $ a_{ij}$ is the $(i,j)$ entry of $A$.
It follows from the relations that $c=h_0+2h_1$ is central. The principal $\mathbb{Z}$-gradation of $\ga$ is given by
\begin{gather*}
\deg h_i = 0,\qquad
\deg e_i = 1,\qquad
\deg f_i = -1.
\end{gather*}
 Let $\gt=\hat{\mathfrak{g}}\oplus \mathbb{C} d$.
\begin{gather*}
\hti_0 = \Span\bigl\{x\in \gt \mid \deg x=0\bigr\},\qquad
\n_\pm = \Span\bigl\{x\in \gt \mid \pm(\deg x)>0\bigr\}.
\end{gather*}
 Then
$
 \gt=\n_-\oplus \hti_0\oplus\n_+$.

A straightforward calculation shows that
\begin{alignat*}{3}
& h_0 \mapsto (-2)\bigl(-3+\sqrt{3}{\rm i}\bigr)X(0) + \frac{2}{3}c,\qquad &&
h_1 \mapsto \bigl(-3+\sqrt{3}{\rm i}\bigr)X(0) + \frac{1}{6}c,&\\
& e_0 \mapsto \frac{4{\rm i}}{\sqrt{3}}\alpha(1)
 + 2\bigl(-1-\sqrt{3}{\rm i}\bigr)X(1),\qquad &&
e_1 \mapsto \frac{-{\rm i}}{\sqrt{3}}\alpha(1)
 + \bigl(-1-\sqrt{3}{\rm i}\bigr)X(1),&\\
& f_0 \mapsto \frac{-2{\rm i}}{\sqrt{3}}\alpha(-1)
 + \bigl(-1-\sqrt{3}{\rm i}\bigr)X(-1),\qquad &&
f_1 \mapsto \frac{{\rm i}}{\sqrt{3}}\alpha(-1)
 + \bigl(-1-\sqrt{3}{\rm i}\bigr)X(-1).&
\end{alignat*}
 extends to a graded Lie algebra isomorphism between $\hat{\mathfrak{g}}$ and $\mathfrak{k}$.
 See also \cite{K90}.

\section{Standard modules}\label{S3}

Let \smash{$\Lambda\in \bigl(\hti_0\bigr)^*$} be a dominant integral weight and let $L(\Lambda)$ be the standard $\mathfrak{k}$-module with highest weight $\Lambda$. Define $\Lambda_0$, $\Lambda_1$ such that $\Lambda_i(h_j)=\delta_{ij}$. The level of a module is $\Lambda(c)=\Lambda(h_0)+2\Lambda(h_1)$. There is only one, up to equivalence, level 1 standard module: $U=L(\Lambda_0)$.
There are
three level 5 standard modules: $L(5\Lambda_0)$, $L(3\Lambda_0+\Lambda_1)$, $L(\Lambda_0+2\Lambda_1)$ and they can be realized inside the tensor product $U^{\otimes 5}$.

Let \smash{$\rho\in \bigl(\hti_0\bigr)^*$} be such that $\rho(h_0)=\rho(h_1)=1$ and $\phi=\Lambda+\rho$. Let
\[
J_{\Lambda}= \{n\in \mathbb{N} \mid n\not\equiv 0,\, \phi(c),\, \pm\phi(h_0),\, \pm\phi(h_1),\, \pm\phi(h_0+h_1) \ (\bmod \, 2\phi(c))  \}
\]
and
\[
K_\Lambda=\begin{cases}
 \{n\in \N \mid n\equiv \pm \phi(h_0) \pmod{2\phi(c)}\}& \text{if }\Lambda(h_0)=\Lambda(h_1),\\
 \varnothing &\text{otherwise}.
\end{cases}
\]
Using Lepowsky's numerator formula and the principal specialization of the Weyl--Kac character formula, we have
\[
\chi_{\Lambda}(q)=\prod_{\underset{n\equiv\pm1 \ ({\rm mod}\, 6)}{n\in \N}}(1-q^n)^{-1}\prod_{n\in J_\Lambda}(1-q^n)^{-1}\prod_{n\in K_\Lambda}(1-q^n).
\]

Set \smash{$H(q)=\prod_{\underset{n\equiv\pm1 \ ({\rm mod}\, 6)}{n\in \N}}(1-q^n)^{-1}$}. For level~5,
we get
\begin{align}
&\chi_{5\Lambda_0}(q)
= H(q)
\prod_{\substack{n\in\N\\
n\equiv \pm 2,\pm 3,\pm 4,\pm 5 \ ({\rm mod}\, {16})}}
(1-q^n)^{-1}, \label{83}\\
&\chi_{3\Lambda_0+\Lambda_1}(q)
= H(q)
\prod_{\substack{n\in\N\\
n\equiv \pm 1,\pm 3,\pm 5,\pm 7 \ ({\rm mod}\, {16})}}
(1-q^n)^{-1},\nonumber\\\ 
&\chi_{\Lambda_0+2\Lambda_1}(q)
= H(q)
\prod_{\substack{n\in\N\\
n\equiv \pm 1,\pm 4,\pm 6,\pm 7\ ({\rm mod}\,{16})}}
(1-q^n)^{-1}. \nonumber 
\end{align}

\begin{Remark}
 These can be compared with formulas (83), (84), and (86), respectively, in \cite{LM78}.
\end{Remark}

We consider the case $\Lambda=5\Lambda_0$.
The product \smash{$\prod_{\substack{n \in \mathbb{N} \\ n \equiv \pm2, \pm3, \pm4, \pm5 \ ({\rm mod}\,{16})}} (1 - q^n)^{-1}$}, when expanded, begins as
\begin{gather*}
1+q^2+q^3+2 q^4+2 q^5+3 q^6+3 q^7+5 q^8+5 q^9+7 q^{10}+8 q^{11}+11 q^{12}+12q^{13}\\
\qquad{}+16 q^{14}+18 q^{15}+23 q^{16}+26 q^{17}+33 q^{18}+37
q^{19}+46 q^{20}+52 q^{21}+63 q^{22}\\
\qquad{}+72 q^{23}+87 q^{24}+98 q^{25}+117 q^{26}+133 q^{27}+157 q^{28}+178 q^{29}+209 q^{30}+236 q^{31}\\
\qquad{}+276 q^{32}+312
q^{33}+361 q^{34}+408 q^{35}+471 q^{36}+530 q^{37}+609 q^{38}\\
\qquad{}+686 q^{39}+784 q^{40}+881 q^{41}+1004 q^{42}+1126 q^{43}+1279 q^{44}+1433 q^{45}\\
\qquad{}+1621
q^{46}+1814 q^{47}+2048 q^{48}+2286 q^{49}+2574 q^{50}+\cdots.
\end{gather*}

\section[Relations on L(5 Lambda\_0)]{Relations on $\boldsymbol{L(5\Lambda_0)}$}\label{S4}

	Define the polynomial
	\begin{equation*}
		P_0(x) = (1-\omega x)\bigl(1-\omega^2x\bigr)^3(1+x)^4\bigl(1-\omega^4x\bigr)^3\bigl(1-\omega^5 x\bigr).
	\end{equation*}
	Let $i_1,i_2\in\mathbb{Z}$ such that $|i_1-i_2| \leq 1$. By \cite[Proposition 3.3]{C92}, the limit on the right-hand side of the following expression exists, and we define
	\begin{equation*}
		X\bigl(\nu^{i_1}\alpha,\nu^{i_2}\alpha;w\bigr) =
		\lim_{w_1,w_2\to w}P_0(w_2/w_1)X\bigl(\nu^{i_1}\alpha;w_1\bigr)X\bigl(\nu^{i_2}\alpha;w_2\bigr).
	\end{equation*}
	Similarly, let $i_1,i_2,i_3\in\mathbb{Z}$ such that $|i_a-i_b| \leq 1$ for all $a$, $b$. Then by \cite[Proposition~3.3]{C92}, the limit on the right-hand side of the following expression exists, and we define
	\begin{align*}
& X\bigl(\nu^{i_1}\alpha,\nu^{i_2}\alpha,\nu^{i_3}\alpha;w\bigr)\\
& \qquad =
\lim_{w_1,w_2,w_3\to w}
\prod_{1\le a<b\le 3} P_0(w_b/w_a)
X\bigl(\nu^{i_1}\alpha;w_1\bigr)
X\bigl(\nu^{i_2}\alpha;w_2\bigr)
X\bigl(\nu^{i_3}\alpha;w_3\bigr).
\end{align*}
	Recall the constant
	\begin{equation*}
		K = \frac{1+\omega}{36},
	\end{equation*}
	from \eqref{vertexOperatorFormula}.
	The following is a special case of \cite[Theorem 4.6]{C92}.
	\begin{Theorem}\label{Rrelation}
		On any level $5$ standard module one has
		\begin{equation*}
			X(\alpha,\alpha,\alpha;w) - 3P_0(1)^2 K E^-(-\alpha;w)X(-\alpha,-\alpha;w)
			E^+(-\alpha;w) = 0.
		\end{equation*}
	\end{Theorem}
	
	The following can be proved by the same method as \cite[Theorem 4.6]{C92}.
	\begin{Theorem}\label{Srelation}
		On any level $5$ standard module one has
		\begin{equation*}
			X(\alpha,\alpha,\nu \alpha; w) - E^-(-\alpha;w)X\bigl(-\alpha,-\alpha,\nu^2 \alpha;w\bigr) E^+(-\alpha;w) = 0.
		\end{equation*}
	\end{Theorem}
	
	Let $\bar {\mathcal{U}}$ be the completion of the enveloping algebra $\mathcal{U}(\mathfrak{k})$ with respect to the family $\mathcal{F}$ of all level 5 Verma modules as in \cite{MP87}. We define the following generating functions and coefficients over respectively in $\bar {\mathcal{U}}$:
	\begin{gather*}
		R(w) = X(\alpha,\alpha,\alpha;w) - 3P_0(1)^2K E^-(-\alpha;w)X(-\alpha,-\alpha;w)
		E^+(-\alpha;w),
	\\
		S(w) = X(\alpha,\alpha,\nu\alpha;w) - E^-(-\alpha;w)X\bigl(-\alpha,-\alpha,\nu^2 \alpha;w\bigr) E^+(-\alpha;w),
	\\
		G(w) = X(\alpha,\alpha,\nu\alpha;w) = \sum_{n\in\mathbb{Z}} G(n)w^{-n},
	\\
		H(w) = X\bigl(-\alpha,-\alpha,\nu^2 \alpha;w\bigr) = \sum_{n\in\mathbb{Z}} H(n)w^{-n},
	\\
		R(w) = \sum_{n\in\mathbb{Z}} R(n)w^{-n},
	\qquad
		S(w) = \sum_{n\in\mathbb{Z}} S(n)w^{-n}.
	\end{gather*}
	
	In order to define ``leading terms'', we introduce a well-order on the Poincar\'e-Birkhoff-Witt-type vectors. We order the basis elements of $\mathfrak{k}$ by $\cdots < \alpha(-7) <\alpha(-5) <\alpha(-1) <\dots < X(-4) <X(-3) <X(-2) <X(-1)< X(0) < \cdots$. For a partition $\pi = (\pi_1,\dots,\pi_j)$, define the \textit{weight} of $\pi$ by
$|\pi| = \sum_{i=1}^j \pi_i$,
	and the \textit{length} of $\pi$ by
$\ell(\pi) = j$.
	For partitions $\lambda = (\lambda_1,\dots,\lambda_h)$, $\lambda_i\equiv \pm1\pmod6$, $\mu = (\mu_1,\dots,\mu_k)$, set
	\begin{gather*}
	\alpha(-\lambda) = \alpha(-\lambda_1)\cdots \alpha(-\lambda_h),
	\qquad
	X(-\mu) = X(-\mu_1)\cdots X(-\mu_k).
	\end{gather*}
	
	Denote by $v_\Lambda$ a fixed highest weight vector in $L(\Lambda)$.
	As the relations $R(-n)$, $S(-n)$ are homogeneous, we shall only need to compare monomials of equal degree. Let
	\begin{equation*}
		\alpha(-\lambda)X(-\mu)v_\Lambda, \qquad \alpha(-\kappa)X(-\pi)v_\Lambda
	\end{equation*}
	be two monomials in $L(\Lambda)$ of degree $-n = -|\lambda|-|\mu| = -|\kappa|-|\pi|$.
	We then define that
	\begin{equation*}
		\alpha(-\lambda)X(-\mu)v_\Lambda < \alpha(-\kappa)X(-\pi)v_\Lambda,
	\end{equation*}
	if one of the following holds:
	\begin{enumerate}\itemsep=0pt
	\item[(1)] $|\mu| > |\pi| $,
	
	\item[(2)] $|\mu| = |\pi|$, and $\ell(\mu) > \ell(\pi)$,
	
	\item[(3)] $|\mu| = |\pi|$, $\ell(\mu) = \ell(\pi)$, and $\mu < \pi$ in lexicographic order,
	
	\item[(4)] $\mu = \pi$, and $\ell(\lambda) > \ell(\kappa)$,
	
	\item[(5)] $\mu = \pi$, $\ell(\lambda) = \ell(\kappa)$, and $\lambda < \kappa$ in lexicographic order.
	\end{enumerate}
	In a linear combination
	\begin{equation*}
		\sum_{i=1}^m c_i\alpha\bigl(-\lambda^{(i)}\bigr)X\bigl(-\mu^{(i)}\bigr)v_\Lambda,
	\end{equation*}
	where \smash{$\bigl|\lambda^{(i)}\bigr| + \bigl|\mu^{(i)}\bigr| = n$}, $c_i\in\mathbb{C}^\times$ for all $i$, we say that \smash{$c_1\alpha\bigl(-\lambda^{(1)}\bigr)X\bigl(-\mu^{(1)}\bigr)v_\Lambda$} is the \textit{leading term} of the sum if{\samepage
 \begin{equation*}
		\alpha\bigl(-\lambda^{(1)}\bigr)X\bigl(-\mu^{(1)}\bigr)v_\Lambda < 	\alpha\bigl(-\lambda^{(i)}\bigr)X\bigl(-\mu^{(i)}\bigr)v_\Lambda
	\end{equation*}
	for all $i > 1$.}
	
	For $n,j \in\mathbb{Z}$, we define
	\begin{equation*}
		L_{(n,j)} = \operatorname{span}\{\alpha(-\lambda)X(-\mu)v_\Lambda \mid |\mu| < n \text { or }\ell(\mu)\leq j\}.
	\end{equation*}
		
	Denote by $\Psi(x)$ the power series expansion into nonnegative integral powers of $x$ of
	\begin{equation*}
		\Psi(x) = \frac{1-3x+4x^2-3x^3+x^4}{1+3x+4x^2+3x^3+x^4} = 1-6x + 18x^2 +\cdots.
	\end{equation*}
	\begin{Proposition}[{\cite[equation (2.2.47]{N14})}]\label{N14}
		We have
		\begin{equation*}
			X(\alpha;w_1)E^-(-\alpha;w_2) = \Psi\biggl(\frac{w_2}{w_1}\biggr)E^-(-\alpha;w_2)X(\alpha;w_1),
		\end{equation*}
		and
		\begin{equation*}
			E^+(-\alpha;w_1)X(\alpha;w_2) = \Psi\biggl(\frac{w_2}{w_1}\biggr) X(\alpha;w_2) E^+(-\alpha;w_1).
		\end{equation*}
	\end{Proposition}
	
	Denote the coefficients in $\Psi(x)$ by $a_i$ so that
	\begin{equation*}
		\Psi(x) = \sum_{i=0}^\infty a_i x^i = 1-6x + 18x^2 +\cdots.
	\end{equation*}
	
	Combining Theorem~\ref{Srelation} with Proposition~\ref{N14}, we obtain the following identity that will be used to obtain leading terms of length~4.
	\begin{Proposition}\label{Sright}
For $p,q\in\mathbb{Z}$,
		\begin{gather}\label{Srighteq}
				S(-p)X(-q)v_\Lambda \equiv G(-p)X(-q) v_\Lambda
				- \sum_{i=0}^\infty a_i H(-p-i)X(-q+i) v_\Lambda \mod L_{(p+q,3)}.
		\end{gather}
	\end{Proposition}
	\begin{proof}
		Using Nandi's identity, Proposition~\ref{N14}, we obtain
		\begin{align*}
 S(w_1)X(\alpha;w_2) &= \bigl(G(w_1) - E^-(-\alpha;w_1)H(w_1)E^+(-\alpha;w_1)\bigr)X(\alpha;w_2) \\
				&= G(w_1)X(\alpha;w_2) - E^-(-\alpha;w_1)H(w_1)X(\alpha;w_2)E^+(-\alpha;w_1)\Psi\biggl(\frac{w_2}{w_1}\biggr).
			\end{align*}
		The coefficient of \smash{$w_1^pw_2^q$} gives \eqref{Srighteq}.
	\end{proof}

	One has the following analogue using multiplication by $X(-q)$ from the left instead of from the~right.
	\begin{Proposition}\label{Sleft}
For $p,q \in\mathbb{Z}$,
		\begin{gather*}
				X(-q)S(-p) v_\Lambda \equiv X(-q)G(-p) v_\Lambda
				- \sum_{i=0}^\infty a_i X(-q-i)H(-p+i) v_\Lambda
				 \mod L_{(p+q,3)}.
		\end{gather*}
	\end{Proposition}
	The proof is similar to that of Proposition~\ref{Sright}.
	More generally, we have the following.
	
	\begin{Proposition}\label{Sleftright}
		For $q_1,\dots,q_h \in\mathbb{Z}$, $r_1,\dots,r_m\in\mathbb{Z}$, $p\in\mathbb{Z}$,
		\begin{gather*}
				X(-q_1)\cdots X(-q_h)S(-p)X(-r_1)\cdots X(-r_m)v_\Lambda \\
				\qquad{}\equiv 	X(-q_1)\cdots X(-q_h)G(-p)X(-r_1)\cdots X(-r_m)v_\Lambda \\
				\qquad\quad{}-\sum_{i_1,\dots,i_h \ge 0}\sum_{j_1,\dots,j_m\ge 0}\prod_{u=1}^h a_{i_u}\prod_{v=1}^m a_{j_v} X(-q_1-i_1)\cdots X(-q_h-i_h) \\
				 \qquad\qquad{}\times H\Bigl(-p+\sum i_u - \sum j_v\Bigr)X(-r_1+j_1)\cdots X(-r_m+j_m)v_\Lambda \mod L_{(s,h+m+2)},
		\end{gather*}
		where $s = p + \sum q_u + \sum r_v$.
	\end{Proposition}
	The proof is analogous to that of Proposition~\ref{Sright}.
	
\section[Calculation of a leading term at degree -14]{Calculation of a leading term at degree $\boldsymbol{-14}$}\label{S5}

	In this section, we exemplify our methods by showing that $X(-(5,5,2,2))v_\Lambda$ is a leading term at degree $-14$. This leads to condition 13 in the list in Section \ref{S6}.
	
	The coefficient of $X(-(a,b,c))$ in $X\bigl(\nu^h\alpha,\nu^i\alpha,\nu^j\alpha;w\bigr)$ (when expanded in the PBW-type monomials) is
	\begin{equation*}
		C_1(h,i,j,a)=\omega^{-a(h+i+j)}
	\end{equation*}
	if $a=b=c$, is
	\begin{equation*}
		C_2(h,i,j,a,b) = \omega^{-a(i+j)-bh} + \omega^{-a(h+j)-bi} + \omega^{-a(h+i)-bj}
	\end{equation*}
	if $a = c$ and $a\neq b$, and is
	\begin{align*}
			C_3(h,i,j,a,b,c)={}& \omega^{-ch-bi-aj} + \omega^{-bh-ci-aj} + \omega^{-ch-ai-bj} \\
			&{+}\, \omega^{-ah-ci-bj} + \omega^{-bh-ai-cj} + \omega^{-ah-bi-cj}
	\end{align*}
if $a$, $b$, $c$ are distinct. Here we have renormalized the $X\bigl(\nu^h\alpha,\nu^i\alpha,\nu^j\alpha;w\bigr)$
 by dividing all coefficients by $P_0(1)^3$. Using these coefficients, we find the following relations at degree $-14$.	
As the coefficient of $w^{14}$ in $X(-\alpha,-\alpha;w)$ in Theorem~\ref{Rrelation} consists of second degree monomials in the~$X(-n)$, we obtain after multiplication of $R(-12)$ by $X(-2)$ the following ``length 4'' terms:
\begin{align}
R(-12)X(-2)v_\Lambda={} (&X(-(4,4,4,2)) + 6 X(-(5,4,3,2)) + 3 X(-(5,5,2,2))\nonumber \\
& +  3 X(-(6,3,3,2)) + 6 X(-(6,4,2,2)) + \cdots ) v_\Lambda.\label{RXcase}
\end{align}
	The case $p = 11$, $q = 3$ of Proposition~\ref{Sright} gives
	\begin{align}
S(-11)X(-3)v_\Lambda = (&{-}4X(-(4,4,3,3)) - 6\omega X(-(4,4,4,2)) \nonumber\\
		& -  3 X(-(5,3,3,3)) + (-2-24\omega) X(-(5,4,3,2))
		 \nonumber\\
		 & +  (6-6\omega) X(-(5,5,2,2)) - 8 X(-(6,3,3,2)) + \cdots )v_\Lambda .\label{SXcase}
	\end{align}
	The case $p = 8$, $q = 6$ of Proposition~\ref{Sleft} gives
	\begin{equation}\label{XScase}
		X(-6)S(-8)v_\Lambda = ((1-2\omega) X(-(6,3,3,2)) + (1-2\omega) X(-(6,4,2,2)) + \cdots ) v_\Lambda.
	\end{equation}
	Here in \eqref{RXcase}--\eqref{XScase} we have left out partitions that have a part $\geq 7$, which are higher in lexicographic order.
	
	In Table~\ref{table1}, we present an $(8\times 7)$-matrix, where we collect the coefficients of the monomials $X(-(a,b,c,d))v_\Lambda$ with all parts $a$, $b$, $c$, $d$ $\leq 6$, in the following expressions:
	\begin{alignat*}{5}
		&\frac13R(-11)X(-3)v_\Lambda,\qquad&& R(-12)X(-2)v_\Lambda,\qquad&& S(-11)X(-3)v_\Lambda,\qquad&& S(-12)X(-2)v_\Lambda,& \\
		&\frac13 X(-6)R(-8)v_\Lambda,\qquad&& X(-5)R(-9)v_\Lambda,\qquad&& X(-6)S(-8)v_\Lambda,\qquad&& X(-5)S(-9)v_\Lambda.&
	\end{alignat*}

	 Note that all partitions that contain a part $1$ can be, and have been, left out, since $f_1.v_\Lambda = 0$ shows that $X(-1)v_\Lambda \in\mathbb{C}\alpha(-1)v_\Lambda$.
	
\begin{table}[!ht]
\centering
\caption{Relations at $n = 14$.} \label{table1}\vspace{1mm}

\footnotesize
\setlength{\tabcolsep}{3pt}
\renewcommand{\arraystretch}{1.05}

\begin{tabular}{*{7}{>{\centering\arraybackslash}p{1.9cm}}}
\toprule
\shortstack{\scriptsize $X(-(4433))v_\Lambda$} &
\shortstack{\scriptsize $X(-(4442))v_\Lambda$} &
\shortstack{\scriptsize $X(-(5333))v_\Lambda$} &
\shortstack{\scriptsize $X(-(5432))v_\Lambda$} &
\shortstack{\scriptsize $X(-(5522))v_\Lambda$} &
\shortstack{\scriptsize $X(-(6332))v_\Lambda$} &
\shortstack{\scriptsize $X(-(6422))v_\Lambda$} \\
\midrule
$\hphantom{-}1$ & $\hphantom{-}0$ & $\hphantom{-}1$ & $\hphantom{-}2$ & $\hphantom{-}0$ & $\hphantom{-}2$ & $0$ \\
$\hphantom{-}0$ & $\hphantom{-}1$ & $\hphantom{-}0$ & $\hphantom{-}6$ & $\hphantom{-}3$ & $\hphantom{-}3$ & $6$ \\
$-4$ & $-6\omega$ & $-3$ & $-2-24\omega$ & $6-6\omega$ & $-8$ & $0$ \\
$\hphantom{-}0$ & $-1+2\omega$ & $\hphantom{-}0$ & $-4+8\omega$ & $-1+2\omega$ & $\hphantom{-}0$ & $0$ \\
$\hphantom{-}0$ & $\hphantom{-}0$ & $\hphantom{-}0$ & $\hphantom{-}0$ & $\hphantom{-}0$ & $\hphantom{-}1$ & $1$ \\
$\hphantom{-}0$ & $\hphantom{-}0$ & $\hphantom{-}1$ & $\hphantom{-}6$ & $\hphantom{-}3$ & $\hphantom{-}0$ & $0$ \\
$\hphantom{-}0$ & $\hphantom{-}0$ & $\hphantom{-}0$ & $\hphantom{-}0$ & $\hphantom{-}0$ & $1-2\omega$ & $1-2\omega$ \\
$\hphantom{-}0$ & $\hphantom{-}0$ & $-2$ & $-8$ & $-1$ & $-18+6\omega$ & $-12+6\omega$ \\
\bottomrule
\end{tabular}
\end{table}

Reducing this matrix to row echelon form yields the matrix in Table~\ref{tab2}.
Thus, the fifth row shows that a linear combination of relations has leading term $X(-(5,5,2,2))v_\Lambda$.
	
\begin{table}[!ht]
\centering
\caption{Row reduced matrix at $n = 14$.}\label{tab2}\vspace{1mm}

\footnotesize
\setlength{\tabcolsep}{3pt}
\renewcommand{\arraystretch}{1.05}

\begin{tabular}{*{7}{>{\centering\arraybackslash}p{1.9cm}}}
\toprule
\shortstack{\scriptsize $X(-(4433))v_\Lambda$} &
\shortstack{\scriptsize $X(-(4442))v_\Lambda$} &
\shortstack{\scriptsize $X(-(5333))v_\Lambda$} &
\shortstack{\scriptsize $X(-(5432))v_\Lambda$} &
\shortstack{\scriptsize $X(-(5522))v_\Lambda$} &
\shortstack{\scriptsize $X(-(6332))v_\Lambda$} &
\shortstack{\scriptsize $X(-(6422))v_\Lambda$} \\
\midrule
$1$ & $0$ & $0$ & $0$ & $0$ & $0$ & $\hphantom{-}4$ \\
$0$ & $1$ & $0$ & $0$ & $0$ & $0$ & $-6$ \\
$0$ & $0$ & $1$ & $0$ & $0$ & $0$ & $-9$ \\
$0$ & $0$ & $0$ & $1$ & $0$ & $0$ & $\hphantom{-}\tfrac{3}{2}$ \\
$0$ & $0$ & $0$ & $0$ & $1$ & $0$ & $\hphantom{-}0$ \\
$0$ & $0$ & $0$ & $0$ & $0$ & $1$ & $\hphantom{-}1$ \\
$0$ & $0$ & $0$ & $0$ & $0$ & $0$ & $\hphantom{-}0$ \\
$0$ & $0$ & $0$ & $0$ & $0$ & $0$ & $\hphantom{-}0$ \\
\bottomrule
\end{tabular}
\end{table}	

	\section{List of leading terms of relations}\label{S6}

	As in \cite{C92,LW82,MP87,N14} and other papers on the structure of standard modules
	over affine Lie algebras, we shall reduce a spanning set of Poincar\'e--Birkhoff--Witt-type vectors by removing the leading terms of relations.
	
	Consider the following conditions on a partition $\mu$ appearing in a PBW-type spanning element~$\alpha(-\lambda)X(-\mu)v_\Lambda \in L(5\Lambda_0)$.
	
	\textbf{List of conditions on $\boldsymbol{\mu}$.}
	\begin{enumerate}\itemsep=0pt
	\item[(1)] no part is equal to 1,
	
	\item[(2)] $\mu$ does not have three equal parts,
\end{enumerate}
in addition $\mu$ does not have a sub-partition of any of the following forms:
	\begin{enumerate}\itemsep=0pt
	\item[(3)] $(k,k,k-1)$,
	
	\item[(4)] $(k,k-1,k-1)$,
	
	\item[(5)] $(k,k-1,k-2)$ with $k \not\equiv 1\pmod6$,
	
	\item[(6)] $(k,k,k-2)$ with $k \not\equiv 4\pmod6$,
	
	\item[(7)] $(k,k-2,k-2)$ with $k \not\equiv 4\pmod6$,
	
	\item[(8)] $(k,k-2,k-3)$ with $k\equiv 5\pmod6$,
	
	\item[(9)] $(k,k-1,k-3)$ with $k \equiv 4\pmod6$,
	
	\item[(10)] $(k,k,k-3)$ with $k \equiv 1 \pmod6$,
	
	\item[(11)] $(4+6k, 4+6k, 2+6k, 2+6k)$,
	
	\item[(12)] $(5+6k, 4+6k, 2+6k, 2+6k)$,
	
	\item[(13)] $(5+k, 5+k, 2+k, 2+k)$,
	
	\item[(14)] $(10+6k, 10+6k, 8+6k, 7+6k)$,
	
	\item[(15)] $(7+6k, 6+6k, 5+6k, 3+6k)$,
	
	\item[(16)] $(7+2k,6+2k, 4+2k, 3+2k)$,
	
	\item[(17)] $(6+k, 5+k, 2+k, 2+k), k\not\equiv 4 \pmod6$,
	
	\item[(18)] $(8+k, 8+k, 5+k, 4+k), k\not\equiv 4 \pmod6$,
	
	\item[(19)] $(8+6k, 8+6k, 5+6k, 3+6k)$,
	
	\item[(20)] $(9+6k, 7+6k, 6+6k, 5+6k)$,
	
	\item[(21)] $(10+6k, 8+6k, 8+6k, 5+6k)$,
	
	\item[(22)] $(7+6k, 4+6k, 4+6k, 2+6k)$,
	
	\item[(23)] $(7+6k, 6+6k, 4+6k, 2+6k, 2+6k)$,
	
	\item[(24)] $(8+6k, 7+6k, 5+6k, 4+6k, 2+6k)$,
	
	\item[(25)] $(9+6k, 9+6k, 5+6k, 5+6k, 2+6k)$,
	
	\item[(26)] $(10+6k,8+6k,7+6k,5+6k,4+6k)$,
	
	\item[(27)] $(10+6k,10+6k,8+6k,6+6k,5+6k)$,
	
	\item[(28)] $(11+k,10+k,8+k,6+k,5+k), k\equiv 0,1,2 \pmod6$,
	
	\item[(29)] $(10+6k,10+6k,8+6k,5+6k,5+6k)$,
	
	\item[(30)] $(11+6k,10+6k,8+6k,5+6k,5+6k)$,
	
	\item[(31)] $(9+6k, 8+6k, 6+6k, 4+6k, 2+6k, 2+6k)$,
	
	\item[(32)] $(10+6k, 10+6k, 8+6k, 6+6k, 4+6k, 3+6k)$,
	
	\item[(33)] $(11+6k, 8+6k, 8+6k, 5+6k, 2+6k, 2+6k)$,
	
	\item[(34)] $(12+6k, 11+6k, 9+6k, 7+6k, 4+6k, 4+6k)$,
 \end{enumerate}
where in all conditions $k\in\mathbb{Z}$.

Using Theorems \ref{Rrelation} and \ref{Srelation} and Propositions \ref{Sright}--\ref{Sleftright}, we have been able to reduce the spanning set of $L(5\Lambda_0)$ to the set of vectors
	$\alpha(-\lambda)X(-\mu)v_\Lambda$ such that $\mu$ satisfies the conditions (1), (2), (3)--(10) with $k \ge 3$ and (11)--(34) with $k=0$, and the sub-partition appearing at the end of $\mu$.
	
In addition, we conjecture that the spanning set can be reduced by removing the vectors $\alpha(-\lambda)X(-\mu)v_\Lambda$, where the partition $\mu$ has a sub-partition, not necessarily at the end of $\mu$, of the form in conditions (3)--(34) with any $k \ge 0$.
	
Condition (1) follows from $f_1. v_\Lambda = 0$.
Conditions (2)--(10) follow by taking linear combinations of $R(-n)v_\Lambda=0$, $S(-n)v_\Lambda=0$ for appropriate $n$.
Conditions (11)--(22) are analogously obtained from the leading terms of relations obtained as linear combinations of the ``length 4'' relations $R(-p)X(-n+p)v_\Lambda$, $S(-p)X(-n+p)v_\Lambda$, $X(-n+p)S(-p)v_\Lambda$. These are obtained by the \textsc{Maple} program ``congruences4.mw''. Conditions (23)--(30) are similarly obtained in ``congruences5.mw''.
Finally, conditions (31)--(34) are obtained in ``congruences6.mw''.
The \textsc{Maple} programs are available for download at
\href{https://arxiv.org/abs/2511.12284}{arXiv:2511.12284}.
	
Using computer programs, we enumerated the partitions $\mu$ satisfying conditions~(1)--(34) and formed the truncated generating series
\begin{equation*}
g(q)=\sum_{|\mu|\le 48} q^{|\mu|}.
\end{equation*}
	Comparing this with the count according to the principally specialized character of $L(5\Lambda_0)$, and dividing by the character of the Heisenberg module, which gives
	\begin{equation}\label{chiq}
		\chi(q) = \prod_{j\equiv \pm2,\pm3,\pm4,\pm5\  ({\rm mod}\, 16)}\bigl(1-q^j\bigr)^{-1},
	\end{equation}
	cf.\ \eqref{83}, we have equality between $g(q)$ and $\chi(q)$ for all powers of $q$ up to 48 except that the coefficients in $g(q)$ exceed that of $\chi(q)$ by 1 for the coefficients of $q^{42}$ and $q^{48}$.
	
	\begin{Remark}
		This count uses additional arguments, for $|\mu| \le 42$, analogous to those in \cite[Sections~4 and~6]{N14}, that allow the subpartition to appear anywhere in $\mu$, what Nandi terms ``forbidden patterns''. We omit these since we feel that they do not contribute to the understanding of the underlying combinatorics.
	\end{Remark}
	
	It appears that our methods of calculation are insufficient to obtain the missing partitions of weights $42$ and $48$.
	
	\section[A basis for a level 2 module over sl(2,C)\^{}~]{A basis for a level 2 module over $\boldsymbol{\mathfrak{sl}(2,\mathbb{C})^\sim}$}\label{S7}

	In \cite[Theorem 2.2.1]{Bo02}, Borcea observed in particular that the principally specialized character of $L(5\Lambda_0)$ is related to the $(2,1)$-specialized character of the module $L(2\Lambda_0)$ over the affine Lie algebra of type \smash{$A_1^{(1)}$} as follows:
	\begin{equation*}
		\operatorname{ch}_q^{(1,1)} \bigl(L_{A_2^{(2)}}(5\Lambda_0)\bigr) = \prod_{i\equiv \pm1 \ ({\rm mod} \,6)} \bigl(1-q^i\bigr)^{-1}\operatorname{ch}_q^{(2,1)}\bigl(L_{A_1^{(1)}}(2\Lambda_0)\bigr),
	\end{equation*}
	where \smash{$\operatorname{ch}_q^{(s,t)}$} denotes the $(s,t)$-specialized grading. Here \smash{$\prod_{i\equiv \pm1\pmod6}\bigl(1-q^i\bigr)^{-1}$} is the character of the Fock space of the Heisenberg Lie algebra spanned by the $\alpha(j)$ and $c$, so that the identity equates $\operatorname{ch}(L(2\Lambda_0))$ with the character of the vacuum space of the Heisenberg Lie algebra.
	
For comparison, we present the PBW-type basis obtained in \cite{MP99} and \cite{FKLMM99} for the module~$L(2\Lambda_0)$ in the $(2,1)$-specialized grading.
	
Let $\mathfrak{a} = \mathfrak{sl}(2,\mathbb{C})$ and form the affine Lie algebra
$\hat{\mathfrak{a}} = \mathfrak{a}\otimes \mathbb{C}\bigl[t,t^{-1}\bigr] \oplus \mathbb{C}c$
	such that
	\begin{equation*}
		[u(m),v(n)] = [u,v](m+n) + m\delta_{m+n,0}\langle u,v\rangle c,
	\end{equation*}
	for $u,v\in\mathfrak{a}$,
	and
$c\in \text{center}(\hat{\mathfrak{a}})$,
	where $u(m):=u\otimes t^m$, and $\langle \cdot,\cdot\rangle$ denotes the invariant symmetric bilinear form on $\mathfrak{a}$ such that $\langle \alpha,\alpha\rangle = 2$ for each root $\alpha$.
	
	Let $x$, $h$, $y$ be the basis \smash{$x = \bigl[\begin{smallmatrix} 0&1 \\ 0&0 \end{smallmatrix}\bigr]$}, $h = \bigl[\begin{smallmatrix} 1&\hphantom{-} 0 \\ 0&-1 \end{smallmatrix}\bigr]$, $y = \bigl[\begin{smallmatrix} 0&0 \\ 1&0 \end{smallmatrix}\bigr]$ of $\mathfrak{a}$.
	Then $e_0 = y(1)$, $e_1=x(0)$, $h_0 = -h(0)+c$, $h_1 = h(0)$, $f_0 = x(-1)$, $f_1 = y(0)$ forms a set of Kac--Moody generators of $\hat{\mathfrak{a}}$. We denote by $\Lambda_i$ for $i=0,1$ the linear forms such that $\Lambda_i(h_j) = \delta_{ij}$.

	In the $(2,1)$-specialized grading, we have the basis of $\hat{\mathfrak{a}}$
	\begin{gather*}
			T(-3k) = h(-k), \qquad
			T(-3k-1) = y(-k), \qquad
			T(-3k-2) = x(-k-1),
		\end{gather*}
	$k\in\mathbb{Z}$, together with $c$. For a partition $\mu = (\mu_1,\dots,\mu_\ell)$, set
	\begin{equation*}
		T(-\mu) = T(-\mu_1)\cdots T(-\mu_\ell).
	\end{equation*}
	
\begin{Theorem}
The set of vectors $\{T(-\mu)v_{2\Lambda_0}\}$ where the $\mu_i > 0$ and $\mu$ satisfies the following conditions $1$--$11$ forms a basis of the $\hat{\mathfrak{a}}$-module $L(2\Lambda_0)$:
\begin{enumerate}\itemsep=0pt
\item[$(1)$] no part is equal to $1$,
\item[$(2)$] $\mu$ does not have three equal parts,
\end{enumerate}	
	in addition $\mu$ does not have a subpartition of any of the following forms:
\begin{enumerate}\itemsep=0pt	
\item[$(3)$] $(k,k,k-1)$,
\item[$(4)$] $(k,k-1,k-1)$,
\item[$(5)$] $(k,k-1,k-2)$ with $k \not\equiv 1\pmod3$,
\item[$(6)$] $(k,k,k-2)$ with $k \not\equiv 1\pmod3$,
\item[$(7)$] $(k,k-2,k-2)$ with $k \not\equiv 1\pmod3$,
\item[$(8)$] $(k,k-2,k-3)$ with $k\equiv 2\pmod3$,
\item[$(9)$] $(k,k-1,k-3)$ with $k \equiv 1\pmod3$,
\item[$(10)$] $(k,k,k-3)$ with $k\not\equiv 0 \pmod3$,
\item[$(11)$] $(k,k-3,k-3)$ with $k\not\equiv 0 \pmod3$,
\end{enumerate}
where in conditions $3$--$11$ $k\in\mathbb{Z}$.
\end{Theorem}
	
	This is the case $k_0=0$, $k_1 = 2$, $k = 2$ of \cite[equation (11.2.8)]{MP99},
where we have exchanged the $2\Lambda_1$-module for the $2\Lambda_0$-module by applying the diagram automorphism of $\hat{\mathfrak{a}}$. The instance of Borcea's duality is that
	\begin{equation*}
		\operatorname{ch}_q^{(2,1)}\bigl(L_{A_1^{(1)}}(2\Lambda_0)\bigr) = \chi(q),
	\end{equation*}
	with $\chi(q)$ as in \eqref{chiq}.
	
\section{Conclusion}

One of the starting points of this project was the duality of Borcea relating representations of~\smash{$A_1^{(1)}$} and \smash{$A_2^{(2)}$}. We wanted to see if there is a correspondence between the bases in terms of partitions, but our results show that these are drastically different.
	
The list of leading terms obtained in Section \ref{S6}
	 is incomplete, and unfortunately we have no idea if the complete list is finite or infinite. As the first leading terms of lengths $4$, $5$, $6$ appear at degrees $12$, $21$, $31$, and the missing leading term at degree 42 might be of length $7$, the increase in degrees follows the sequence $9$, $10$, $11$ (?) and might proceed regularly to infinity.
	
	The calculations of linear combinations using computer gets difficult with increasing length as the number of partitions to consider grows fast, and in addition the relations have coefficients in the algebraic number field $\mathbb{Q}(\omega)$, which gives us difficulties to handle with \textsc{Mathematica}.
	
	For $n=42$, we suspect that the missing leading term is $X(-(10,10,8,6,4,2,2))$ of length 7. This is based on the coincidence of bivariate generating functions (taking into account the length of partitions) between our module $L(5\Lambda_0)$ and that of the module over the Virasoro algebra at $(c,h) = (1/2,0)$ as in \cite{AEH20}. This indicates that the missing partition is of length~7, and there is only one such candidate.
	
\appendix 		
	\section{Appendix}

We shall now give some explanations on the \textsc{Maple} programs ``congruences4.mw'', ``congruences5.mw'' and ``congruences6.mw'', that were used in Section~\ref{S6}, and can be downloaded at~\href{http://arxiv.org/abs/2511.12284}{arXiv:2511.12284}.
	In the \textsc{Maple} programs, our main job is to implement Propositions~\ref{Sright}--\ref{Sleftright}.
 In the \textsc{Maple} program ``congruences4.mw'', the procedure
$SX(p,q,{\rm maximal\_part})$
 calculates the terms in  $S(-p)X(-q)v_\Lambda$
 with partitions $[a_1,a_2,a_3,a_4]$ in
$
 	X(-(a_1,a_2,a_3,a_4))v_\Lambda$
 satisfying
\[
{\rm maximal\_part} \ge a_1 \ge a_2 \ge a_3 \ge a_4.
\]
 The relation is presented as a list of terms of the form $[a+b\omega, [a_1,a_2,a_3,a_4]]$, where $a,b\in\mathbb{Z}$, and we have for simplicity omitted the global minus sign. The relation can be printed by the command ${\rm normal\_print}(SX(p,q,{\rm maximal\_part}))$.
 To calculate a set of relations for a given degree, we first assign $n:= {\rm degree}$, and select an appropriate cut-off ``maximalPart''. With the script ``setPartitionsResults'', a set of parameter values is selected that by our experience covers all relations that might contribute. Their coefficients are then collected in a matrix, the analogue of the matrix in Table~\ref{table1}, and this is then reduced to row echelon form. From the pivots, one can then read off the leading terms produced.
	
\subsection*{Acknowledgements}

We thank S.~Tsuchioka for helpful correspondence that led to a correction of the original version of the paper.
The authors are also grateful to the anonymous referees for their careful reading of the manuscript and for several helpful comments and suggestions that improved the paper.
 M.~Primc is partially supported by the Croatian Science Foundation under the project IP-2022-10-9006 and by the project ``Implementation of cutting-edge research and its application as part of the Scientific Center of Excellence for Quantum and Complex Systems, and Representations of Lie Algebras'', grant no. PK.1.1.10.0004, co-financed by the European Union through the European Regional Development Fund -- Competitiveness and Cohesion Programme 2021-2027.


\pdfbookmark[1]{References}{ref}
\LastPageEnding

\end{document}